\documentclass[wssquare]{ppnproc}

\newcommand{\eps}{{\eps}}

\newcommand{\R}{{\mathbb R}}

\newcommand{\Sp}{\mathbb S}

\def\1{\mathbb I}
\def\D{\mathsf D\kern 0.5pt}
\renewcommand{\(}{\left(}
\renewcommand{\)}{\right)}
\renewcommand{\eps}{\varepsilon}

\usepackage{color}
\usepackage[colorlinks=true]{hyperref}
\hypersetup{urlcolor=blue, citecolor=red, linkcolor=blue}
\newcommand{\Email}[1]{E-mail: \href{mailto:#1}{\textsf{#1}}}
\newcommand{\dx}{\,dx}
\newcommand{\BFS}{\mathcal R_{\rm FS}}
\newcommand{\nrmcnd}[2]{\|{#1}\|_{\mathrm L^{#2}(\R\times\Sp^{d-1})}}
\newcommand{\irdmu}[1]{\int_{\R_+\times\Sp^{d-1}}{#1}\,d\mu}

\newcommand*\subjclass[2][2010]{\def\@subjclass{#2}{\vspace*{6pt}\subjclassname\csname subjclassname@#1\endcsname}\scriptsize #2\normalsize}
\newcommand{\subjclassname}{\vspace*{6pt}\hspace*{3pt}\scriptsize{2010} Mathematics Subject Classification:}

\begin{document}
\title{Symmetry of optimizers \\
of the Caffarelli-Kohn-Nirenberg inequalities}

\author{J. Dolbeault$^1$ and M.J. Esteban$^2$}

\address{CEREMADE (CNRS UMR n$^\circ$ 7534),\\ PSL research university, Universit\'e Paris-Dauphine,\\
Place de Lattre de Tassigny, 75775 Paris 16, France\\
$^1$\Email{dolbeaul@ceremade.dauphine.fr}\\
$^2$\Email{esteban@ceremade.dauphine.fr}
}

\author{M. Loss}

\address{School of Mathematics, Georgia Tech,\\
686 Cherry St., Atlanta, GA, 30332-0160, USA\\
\Email{loss@math.gatech.edu}}

\begin{abstract}
In their simplest form, the Caffarelli-Kohn-Nirenberg inequalities are a two parameter family
of inequalities. It has been known that there is a region in parameter
space where the optimizers for the inequalities have broken symmetry. It has been shown
recently that in the complement of this region the optimizers are radially
symmetric. The ideas for the proof will be given.
\end{abstract}



\bodymatter
\thispagestyle{empty}

\section{Introduction}

Symmetries of optimizers in variational problems is a central theme in the calculus of variations.
Sophisticated methods like rearrangement inequalities, reflection methods and moving plane methods belong now to the standard repertoire of any analyst. There are, however, examples where these
methods cannot be applied. Variational problems that depend on parameters very often cannot be treated by such methods, simply because, depending on the parameters, the optimizers
are symmetric and sometimes not. Famous examples are the minimizers of the Ginzburg-Landau functional in superconductivity, where, depending on the strength of the quartic interaction the minimizers form a single, symmetric vortex or a vortex lattice. Clearly such problems
cannot be treated by general methods. For certain parameters they ought to work while in others they cannot. Thus, rather special techniques, tailored to the problems at hand, have to be developed to prove symmetry in the desired regions.

One class of such examples is given by the Caffarelli-Kohn-Nirenberg inequalities~\cite{Caffarelli-Kohn-Nirenberg-84}. In these notes, we shall specifically consider the case of the inequality
\begin{equation}\label{CKN}\tag{CKN}
\,\int_{\R^d}\frac{|\nabla w|^2}{|x|^{2\kern 0.5pt a}}\dx \ge \mathsf C^{\kern 0.25pt d}_{a,b} \(\int_{\R^d}\frac{|w|^p}{|x|^{b\kern 0.5pt p}}\dx\)^{2/p}
\end{equation}
with $a\leq b\leq a+1$ if $d\ge 3\,$, $a< b\leq a+1$ if $d=2\,$, and $a < a_c$ where
\[
a_c:=\frac{d-2}2 \,,\quad p=\frac{2\,d}{d-2+2\,(b-a)} \,.
\]
The function $w$ is in a suitable function space which contains, for instance, all smooth functions with compact support. The constant $\mathsf C^{\kern 0.25pt d}_{a,b}$ is, by definition, the best possible constant. Rotating the function $w$ does not change the value
of the various expressions in~\eqref{CKN}, i.e., the inequality is
rotationally invariant. The special case where $a \ge 0$ has been treated by various authors (see the references in~\cite{DEL}). Rearrangement inequalities can be used to reduce the problem to the set of radial functions, for which the optimality issue can then be solved explicitly.

For the case where $a<0$ the problem is much more subtle. Nevertheless,
Catrina and Wang~\cite{Catrina-Wang-01}, proved that the optimizers, i.e., the functions that yield equality in~\eqref{CKN}, exist in the open strip $a <b< a+1$. This result establishes the existence
of non-negative solutions $w \in\mathrm L^p(\R^d; |x|^{-\kern 0.5pt b\kern 0.5pt p}\dx)$ of the equation
\begin{equation} \label{equation}
-\mathrm{div}\big(\kern 0.5pt|x|^{-2a}\,\nabla w\big) = |x|^{-\kern 0.5pt b\kern 0.5pt p}\,w^{p-1} \,.
\end{equation}
Moreover, in the same paper Catrina and Wang also showed that, in some region in the $(a,b)$ plane, the rotational symmetry of the optimizers is broken.
A more detailed analysis by Felli and Schneider~\cite{Felli-Schneider-03} shows that the region
where the optimizers have a broken symmetry contains the set $\BFS:= \{ (a,b)\,:\,a<0\,,\;b < b_{\kern 0.5pt\rm FS}(a)\}$ where
$$
b_{\kern 0.5pt\rm FS}(a):=\frac{d\,(a_c-a)}{2\,\sqrt{(a_c-a)^2+d-1}}+a-a_c\,.
$$
We call this region $\BFS$ the \emph{Felli-Schneider region}.

In~\cite{Felli-Schneider-03} more is shown. The optimizers in the {\it radial class}
can be determined explicitly which allows to compute the second variation operator about these solutions.
The lowest eigenvalue of this operator is strictly negative for $(a,b) \in \BFS$,
equals zero on the curve $b=b_{\kern 0.5pt\rm FS}(a)$ and is strictly positive in the open complement of the Felli-Schneider region: there,
the radial optimizers are stable. Needless to say
that positivity of the second variation does not imply the radial symmetry of the (global) optimizers for the~\eqref{CKN} inequality. Thus, it is a natural question whether or not the optimizers
possess rotational symmetry in the complement of $\BFS$. Let $2^*:=\frac{2\,d}{d-2}$ if $d\ge3$ and $2^*:=\infty$ if $d=2$. The following
theorem is proved in~\cite{DEL}:

\begin{theorem} \label{CKNequation} Let $d \ge 2$, $p \in (2, 2^*)$, $a<0$ and $b$ in the complement of the Felli-Schneider region and such that $p=\frac{2\,d}{d-2+2\,(b-a)}$.
Then any non-negative solution $w \in\mathrm L^p(\R^d; |x|^{-\kern 0.5pt b\kern 0.5pt p}\dx)$ of~\eqref{equation} must be of the form
$$
\left(A+B\,|x|^{2\alpha}\right)^{-\frac{n-2}{2}}
$$
where $A$, $B$ are positive constants,
\begin{equation} \label{alpha}
\alpha = \frac{(1+a-b)\,(a_c -a)}{a_c-a +b}
\end{equation}
and
\begin{equation} \label{en}
n = \frac{2\,p}{p-2} \,.
\end{equation}
In particular this holds for the optimizers of~\eqref{CKN}.
\end{theorem}
There are some interesting consequences. Using the change of variables
$$
w(r, \omega) = r^{a-a_c}\,\phi\big(\log r, \omega\big) \,,
$$
equation~\eqref{equation} can be cast in the form
\begin{equation} \label{equationline}
-\,\partial_z^2 \phi-\,\Delta_\omega\kern 0.5pt \phi + \Lambda\,\phi = \phi^{p-1}\,.
\end{equation}
Here, $\frac x{|x|}=\omega \in \Sp^{d-1}$, $r=|x|$, $z = \log r$, $\Delta_\omega\kern 0.5pt$ is the Laplace-Beltrami
operator on the sphere $\Sp^{d-1}$ and
$$
\Lambda = (a-a_c)^2 \,.
$$
Thus, $\phi$ is a function on the cylinder $\R \times\Sp^{d-1}$. Moreover, as noticed in~\cite{Catrina-Wang-01},~\eqref{CKN} is transformed into
\begin{equation} \label{CKNcylinder}
\nrmcnd{\partial_z\phi}2^2+\nrmcnd{\nabla_\omega\phi}2^2+\Lambda\,\nrmcnd{\phi}2^2\ge\mathsf C^{\kern 0.25pt d}_{a,b}\,\nrmcnd{\phi}p^2\,.
\end{equation}
\begin{corollary}\label{cylinder}
Let $d \ge 2$, $p \in (2,2^*)$. Any non-negative solution $\phi \in\mathrm L^p(\R \times \Sp^{d-1}; dz \,d\omega)$ of~\eqref{equationline} is, up to translations, of the form
$$
\phi_\Lambda (z)=\left( \tfrac2{p\,\Lambda}{\,\cosh^2\big(\tfrac{p-2}2\,\sqrt\Lambda \, z\big)}\right)^{-\frac1{p-2}}\,,
$$
if and only if
$$
\Lambda \le 4\,\frac{d-1}{p^2-4}\,.
$$
In this range, equality in~\eqref{CKNcylinder} is achieved if and only if $\phi(z)=\phi_\Lambda(z+z_0)$ for some $z_0\in\R$.
\end{corollary}
To put this result in perspective we compare it with a result in~\cite{BV-V}.
\begin{theorem} Let $d \ge 2$, $p \in (2, 2^*)$. 
On the sphere $\Sp^d$ consider the equation
$$
Ð\,\Delta\kern 0.5pt u + \lambda\,u = u^{p-1}
$$
with $\lambda>0$. Here $\Delta$ represents the Laplace-Beltrami operator on $\Sp^d$. Then the constant function $u\equiv\lambda^{1/(p-2)}$ is the only non-negative solution
if and only if
$$
\lambda \le \frac{d}{p-2} \,.
$$
\end{theorem}
Thus, Corollary~\ref{cylinder} can be viewed as an extension of the above mentioned \emph{rigidity} result to the non-compact case of a cylinder. As a special case, this also allows to identify the equality case in the interpolation inequality~\eqref{CKNcylinder} on the cylinder.

In the next sections some ideas about the proof are given: we start by the simple case of the standard Sobolev inequality in Section~\ref{Sec:Sobolev}, explain in Section~\ref{Sec:CKN} how to recast~\eqref{CKN} as a Sobolev type inequality in an artificial \emph{dimension} $n$, where $n$ is not necessarily an integer, and conclude by explaining how the main estimates can be produced using a \emph{fast diffusion} flow.

\section{Heuristics for the proof of Theorem~\ref{CKNequation}}\label{Sec:Sobolev}

In order to avoid long computations it is best to explain the ideas in a `simple' example.
For any $d\ge3$, the Sobolev inequality
\begin{equation} \label{sobolev}
\int_{\R^d} |\nabla u|^2\dx \ge \mathsf C_d \left(\int_{\R^d} |u|^p\dx\right)^{2/p} \,, \quad\mbox{with}\quad p =2^*= \frac{2\,d}{d-2}
\end{equation}
is extremely well understood~\cite{Talenti-76,Aubin-76,Lieb-83}. Once more $\mathsf C_d$ denotes
the sharp constant. Note that this inequality appears as a special case of~\eqref{CKN}
if one sets $a=b=0$, in which case $\mathsf C_d=\mathsf C^{\kern 0.25pt d}_{0,0}$. There is equality in~\eqref{sobolev} if and only if $u$ is a translate
of the Aubin-Talenti function
\begin{equation*} \label{sobolevoptimizer}
\left({c_\star\,\lambda + \frac{|x|^2}{\lambda}}\right)^{-(d-2)/2} \,,
\end{equation*}
where $c_\star$ and $\lambda$ are positive constants.
There have been some proofs using flow methods to understand this inequality~\cite{MR1777035,CL}. The flow used for the case at hand is a porous medium / fast diffusion flow. It is given by
\begin{equation} \label{fastdiffusion}
\frac{\partial\kern 0.5pt v}{\partial t} = \Delta\kern 0.5pt v^{1-\frac1d}
\end{equation}
and has the self-similar solutions
$$
v_\star(x,t) = \left({c_\star\,t + \frac{|x|^2}{t}}\right)^{-d} \,.
$$
This function has slow decay in the $x$ variable.
The obvious similarity of the expressions of the Aubin-Talenti and self-similar functions suggests a reformulation of the Sobolev functional by
setting
$$
v = u^{\frac{2\kern 0.5ptd}{d-2}}\,.
$$
Let us define a \emph{pressure variable} $\mathsf p$ by
$$
v = \mathsf p ^{-d} \,.
$$
A short computation shows
\begin{lemma} The Sobolev inequality, written in terms of $v$ and $\mathsf p$, is given by
\begin{equation} \label{sobolevmodified}
a_c^2\int_{\R^d} v\,|\nabla \mathsf p|^2\dx \ge \mathsf C_d \left(\int_{\R^d} v \dx\right)^\frac{d-2}d\,.
\end{equation}
\end{lemma}
Assume now that $v$ satisfies the fast diffusion equation~\eqref{fastdiffusion}.
This implies that $\mathsf p$ evolves by the equation
$$
\frac{\partial\kern 0.5pt\mathsf p}{\partial t} = \frac{d-1}{d} \left(\mathsf p\,\Delta\kern 0.5pt \mathsf p -d\,|\nabla \mathsf p|^2\right) \,.
$$
The right side of~\eqref{sobolevmodified} does not change if $v$ evolves via
\eqref{fastdiffusion}. For the left side we have
\begin{lemma} Assume that $v$ evolves via~\eqref{fastdiffusion}. Then
\begin{eqnarray*}
\frac{d}{dt} \int_{\R^d} v\,|\nabla \mathsf p|^2\dx &=&-\,2 \int_{\R^d} \left[ \tfrac{1}{2} \,\Delta\kern 0.5pt |\nabla \mathsf p|^2 - \nabla \mathsf p \cdot \nabla \Delta\kern 0.5pt \mathsf p - \tfrac{1}{d}\, (\Delta\kern 0.5pt \mathsf p)^2\right] \mathsf p^{1-d}\dx \\
&=&-\,2 \int_{\R^d}{\rm Tr} \left[ \mathrm H_{\mathsf p } - \tfrac{1}{d}\, ({\rm Tr}\,\mathrm H_{\mathsf p })\,{\rm Id}\right]^2 \mathsf p^{1-d}\dx
\end{eqnarray*}
where $\mathrm H_{\mathsf p }=(\nabla\otimes\nabla)\,\mathsf p$ denotes the Hessian matrix of $\mathsf p$. Moreover,
$$
\mathrm H_{\mathsf p } - \tfrac{1}{d}\, ({\rm Tr}\,\mathrm H_{\mathsf p })\,{\rm Id} = 0
$$ if and only if $\mathsf p(x) = a +b\cdot x +c\,|x|^2$ for some $(a,b,c)\in\R\times\R^d\times\R$.
\end{lemma}
The proof is a somewhat longish but straightforward computation. Note, that it is precisely the particular choice of $v$ and $\mathsf p$ that renders the time derivative in such a simple form.

To summarize, while the right side of the Sobolev inequality stays fixed the left side diminishes under the flow.
{\it The idea is to use the fast diffusion flow to drive the functional towards its optimal value.} Actually we use the fact that if $v$ is optimal in~\eqref{sobolevmodified}, or if it is a critical point, the functional has to be stationary under the action of the flow, which allows to identify $\mathsf p$, hence $v$.
To exploit this idea for the~\eqref{CKN} inequality we have to rewrite it in the form of
a Sobolev type inequality.

\section{A modified Sobolev inequality}\label{Sec:CKN}

The first step in the proof is to rewrite the problem in a form that resembles the Sobolev inequality.
If we write
\[
w(r, \omega) =u(s, \omega)\quad\mbox{with}\quad s=r^\alpha\,,
\]
the inequality~\eqref{CKN} takes the form
\begin{equation*}
\int_{\R_+ \times \Sp^{d-1}} \left[\alpha^2 \left(\frac{\partial u}{\partial s}\right)^2 + \frac{|\nabla_\omega u|^2}{s^2}\right]s^{n-1} \,ds\,d\omega
\ge \mathsf C^{\kern 0.25pt d}_{a,b}\,\alpha^{1-\frac{2}{p}} \left(\int_{\R_+ \times \Sp^{d-1}} |u|^p \,s^{n-1} \,ds\,d\omega \right)^{\frac{2}{p}}
\end{equation*}
where $d\omega$ denotes the uniform measure on the sphere $\Sp^{d-1}$, $\nabla_\omega$ denotes the gradient
on $\Sp^{d-1}$ and where $\alpha$ and $n$ are given by~\eqref{alpha} and~\eqref{en}. We shall
abbreviate
$$
\textstyle\D u:=\left(\alpha\,\frac{\partial u}{\partial s},\frac1s\,\nabla_\omega u\right)\,,\quad |\D u|^2 =\alpha^2 \left(\frac{\partial u}{\partial s}\right)^2 + \frac{|\nabla_\omega u|^2}{s^2} \,.
$$
Our inequality is therefore equivalent to a Sobolev type inequality and takes the form
\begin{equation}\label{sob}
\irdmu{|\D u|^2}\ge \mathsf C^{\kern 0.25pt d}_{a,b}\,\alpha^{1-\frac{2}{p}} \left(\irdmu{|u|^p}\right)^{\frac{2}{p}}\,, \quad\mbox{with}\quad p = \frac{2\,n}{n-2}\,.
\end{equation}
This inequality generalizes~\eqref{sobolev}. Here the measure $d\mu$ is defined on $\R^+\times\Sp^{d-1}$ by
$$
d\mu=s^{n-1} \,ds\,d\omega\,.
$$
As in Section~\ref{Sec:Sobolev}, we may consider $v=u^p$ and define a pressure variable $\mathsf p$ such that $v=\mathsf p^{-n}$, so that $u = \mathsf p^{-(n-2)/2}$. With these notations,~\eqref{sob} can be rewritten as
\begin{equation}\label{sob2}
\tfrac14\,(n-2)^2\irdmu{v\,|\D \mathsf p|^2}\ge \mathsf C^{\kern 0.25pt d}_{a,b}\,\alpha^{1-\frac{2}{p}} \left(\irdmu{v}\right)^{\frac{2}{p}}\,.
\end{equation}
With straightforward abuses of notations, we shall write $\irdmu f=\int_{\R^d}f\,d\mu$ and identify $\mathrm L^p(\R_+ \times \Sp^{d-1};d\mu)$ with $\mathrm L^p(\R^d;|x|^{n-d}\,dx)$ or simply $\mathrm L^p(\R^d;d\mu)$. 

One should note that $n$ is, in general, not an integer and the above inequality reduces to Sobolev's
inequality only if $n=d$. Of particular significance is that the curve
\[
b=b_{\kern 0.5pt\rm FS}(a) \,,
\]
when represented in the new variables $\alpha$ and $n$, is given by the equation $\alpha=\alpha_{\rm FS}$ with
\[
\alpha_{\rm FS}:=\sqrt{\frac{d-1}{n-1}}\,.
\]
Thus, for $\alpha > \alpha_{\rm FS}$ the minimizers are not radial.
The equation~\eqref{equation} transforms into the equation
\begin{equation} \label{sobequn}
-\,\mathcal L\,\, u = u^{p-1} \,,
\end{equation}
where $\mathcal L$ is the Laplacian associated with the quadratic form given by the left side of~\eqref{sob}, i.e., $\mathcal L = -\,\mathsf D^*\cdot \mathsf D$.
Theorem~\ref{CKNequation} can be reformulated as
\begin{theorem} \label{CKNreform} Let $d \ge 2$, $p \in (2, 2^*)$, $n=\frac{2\,p}{p-2}>d$ and $\alpha\le\alpha_{\rm FS}$.
Then any non-negative solution $u \in\mathrm L^p(\R^d;d\mu)$ of~\eqref{sobequn} must be of the form
\begin{equation} \label{talentiform}
\left({A+B\,|x|^2}\right)^{-\frac{n-2}{2}}
\end{equation}
where $A$, $B$ are positive constants, and $n$ is given by~\eqref{en}. As a special case, equality in~\eqref{sob2} is achieved if and only if $u$ is given by~\eqref{talentiform}.
\end{theorem}

The upshot of this work can be summarized in the following fashion: {\it Any optimizer in the radial class that is not
unstable under small perturbations is in fact a global minimizer for the~\eqref{CKN} inequality.}

\section{The flow}\label{Sec:Flow}
We consider the fast diffusion flow
\begin{equation} \label{flow}
\frac{\partial\kern 0.5pt v}{\partial t} = \mathcal L\, v^{1-\frac1n} \,.
\end{equation}
It is easily seen that the flow~\eqref{flow} has the self-similar solutions
$$
v_\star(t;s,\omega) = t^{-n} \left( c_\star + \frac{s^2}{2\,(n-1)\,\alpha^2\,t^2} \right)^{-n} \,.
$$

The basic idea is now quite simple. We consider a non-negative solution $u \in\mathrm L^p(\R^d;d\mu)$ of
\eqref{sobequn} and set $v = u^p$. We also consider the pressure variable $\mathsf p$ such that $v=\mathsf p^{-n}$. The first thing to note is that the right side of~\eqref{sob2} does not change if we evolve~$v$ and hence $u$ under the flow~\eqref{flow}. Further, if we differentiate the left
side of~\eqref{sob2} along the flow we obtain
$$
\frac{d}{dt}\irdmu{v\,|\D \mathsf p|^2}=-\,2\irdmu{\left[\tfrac{1}{2} \,\mathcal L\, |\D \mathsf p|^2 - \D \mathsf p \cdot \D \mathcal L\, \mathsf p - \tfrac{1}{n}\, (\mathcal L\, \mathsf p)^2 \right]}\,.
$$
On the other hand simple computations show
that
\begin{equation} \label{zero}
\tfrac14\,(n-2)^2\,\frac{d}{dt}\!\left(\irdmu{\kern-18ptv\,|\D \mathsf p|^2}\right)\!\mbox{\raisebox{-5pt}{$\Big |_{t=0}$}}=-\,2\irdmu{\kern-18pt(\mathcal L\, u)\, u^{1-p}\,\big(\mathcal L\, u^{p\,(n-1)/n}\big)}
\end{equation}
{\it when expressed in terms of $u$}. Now we take $v=u^p$, where $u$ is the solution to~\eqref{sobequn}, as initial datum for~\eqref{equation}. With this choice, the right side in~\eqref{zero} is actually zero. Indeed, by multiplying both sides of~\eqref{sobequn} by $u^{1-p}\,\big(\mathcal L\, u^{p\,(n-1)/n}\big)$
one obtains
$$
\irdmu{( \mathcal L\, u)\,u^{1-p}\,\big(\mathcal L\, u^{p\,(n-1)/n}\big)}=\irdmu{ u^{p-1}\,u^{1-p}\,\big(\mathcal L\, u^{p\,(n-1)/n}\big)}= 0\,.
$$
The interesting point, and the heart of the argument, is that
$$
0=\irdmu{\left[\tfrac{1}{2} \,\mathcal L\, |\D \mathsf p|^2 - \D \mathsf p \cdot \D \mathcal L\, \mathsf p - \tfrac{1}{n}\, (\mathcal L\, \mathsf p)^2 \right]}
$$
can be written as a sum of non-negative terms
precisely when $\alpha \le \alpha_{\rm FS}$, and the vanishing of these terms shows that $u$ must be of the form
$(A+B\,s^2)^{-(n-2)/2}$. In this way one obtains a classification of the non-negative solutions of~\eqref{sobequn}
provided they are in $\mathrm L^p(\R^d;d\mu)$. To simplify notations, we shall omit the index $\omega$, so that from now on $\nabla$ and $ \Delta$ respectively refer to the gradient and to the Laplace-Beltrami operator on $\Sp^{d-1}$. With the notation $'=\partial_s$, our identity can be reworked as follows.
\begin{lemma} \label{first} 
Assume that $d\ge3$, $n>d$ and let $\mathsf p$ be a positive function in $C^3(\Sp^{d-1})$. Then
\begin{multline*}
\tfrac{1}{2} \,\mathcal L\, |\D \mathsf p|^2 - \D \mathsf p \cdot \D \mathcal L\, \mathsf p - \tfrac{1}{n}\, (\mathcal L\, \mathsf p)^2\\
=\alpha^4\,\tfrac{n-1}{n} \left[ \mathsf p'' - \frac{ \mathsf p'}{r} - \frac{\Delta\kern 0.5pt \mathsf p}{\alpha^2\,(n-1)\,r^2}\right]^2
+ \frac{2\,\alpha^2}{r^2} \Big | \nabla \mathsf p' - \frac{\nabla \mathsf p}{r} \Big |^2\\
+\,\frac{1}{r^4} \left[\tfrac{1}{2} \,\Delta\kern 0.5pt|\nabla \mathsf p|^2 - \nabla \mathsf p \cdot \nabla \Delta\kern 0.5pt \mathsf p - \tfrac{1}{n-1}(\Delta\kern 0.5pt \mathsf p)^2 - (n-2)\,\alpha^2\,|\nabla \mathsf p|^2\right] \,.
\end{multline*}
\end{lemma}
The only term in Lemma~\ref{first} that does not have a sign is the last one. When integrated
against $\mathsf p^{1-n}$ over $\Sp^{d-1}$, however, this term can be written as a sum of
squares. The following lemma holds for $d\ge 3$. For the case $d=2$ we refer the reader
to~\cite{DEL}.
\begin{lemma} \label{second}
Assume that $d\ge3$ and that $\mathsf p$ is a positive function in $C^3(\Sp^{d-1})$. Then
\begin{eqnarray*}
&&\kern-24pt\int_{\Sp^{d-1}} \left[\tfrac{1}{2} \,\Delta\kern 0.5pt|\nabla \mathsf p|^2 - \nabla \mathsf p \cdot \nabla \Delta\kern 0.5pt \mathsf p - \tfrac{1}{n-1}(\Delta\kern 0.5pt \mathsf p)^2 - (n-2)\,\alpha^2\,|\nabla \mathsf p|^2\right] \mathsf p^{1-n} \,d\omega\\
&=&\tfrac{(n-2)\,(d-1)}{(n-1)\,(d-2)}\int_{\Sp^{d-1}} {\left\Vert\mathsf L\kern 0.5pt\mathsf p -\tfrac{3\,(n-1)\,(n-d)}{2\,(n-2)\,(d+1)}\,\mathsf M\kern 0.5pt\mathsf p\right\Vert^2 \mathsf p^{1-n}}\,d\omega\\
&&+\,\tfrac{n-d}{2\,(d+1)} \left[\tfrac{n+3}{2} + \tfrac{3\,(n-1)\,(n+1)\,(d-2)}{2\,(n-2)\,(d+1)} \right] \int_{\Sp^{d-1}} { \frac{|\nabla \mathsf p|^4}{\mathsf p^2} \mathsf p^{1-n}}\,d\omega\\
&&+\,(n-2) \left[\alpha_{\rm FS}^2 - \alpha^2\right] \int_{\Sp^{d-1}}{ |\nabla \mathsf p|^2\,\mathsf p^{1-n}}\,d\omega
\end{eqnarray*}
where $\mathsf L\kern 0.5pt\mathsf p :=(\nabla\otimes\nabla)\,\mathsf p-\frac1{d-1}\,(\Delta\kern 0.5pt \mathsf p)\, g$ and $\mathrm M \kern0.5pt \mathsf p:=\frac{\nabla \mathsf p \otimes\nabla \mathsf p} {\mathsf p} -\frac1{d-1}\,\frac{|\nabla \mathsf p |^2}{ \mathsf p} \,g$. Here $g$ is the standard metric on $\Sp^{d-1}$ and $\mathrm L\kern 0.5pt\mathsf p$ denotes the trace free Hessian of $\mathsf p$.
\end{lemma}
The key device used for the proof of this lemma is the Bochner-Lichnerowicz-Weitzenb\"ock
formula. If $\mathcal M$ is a compact Riemannian manifold, then for any smooth function
$f: \mathcal M \to \R$ we have
$$
\tfrac{1}{2} \,\Delta\kern 0.5pt |\nabla f|^2 = \Vert \mathrm H_f \Vert^2 +{\rm Ric}(\nabla f, \nabla f)
$$
where $\Vert \mathrm H_f \Vert^2$ is the trace of the square of the Hessian of $f$ and
${\rm Ric}(\nabla f, \nabla f)$ is the Ricci curvature tensor contracted against
$\nabla f \otimes \nabla f$. If $\mathcal M =\Sp^{d-1}$, then ${\rm Ric}(\nabla f, \nabla f)
= (d-2)\,|\nabla f|^2 $. The main point in Lemma~\ref{second} is that, provided $\alpha\le\alpha_{\rm FS}$, all terms are non-negative.

It is quite easy to see that the vanishing of these terms entails that $\mathsf p$ can only
depend on the variable $s=|x|$ and must be of the form~\eqref{talentiform}.

While the formal computations are straightforward there is the perennial issue of the boundary terms that occur in all the integration by parts. This is due to the fact that one is dealing with solutions of~\eqref{sobequn} and it is not at all clear that the boundary terms vanish. This requires a detailed regularity analysis
of the solutions of~\eqref{sobequn}. The task is non-trivial because the exponent $p$ is
critical for the scaling in the $s$ variable. The reader may consult~\cite{DEL} for details.

The computations outlined above can be carried over to the case where $\Sp^{d-1}$
is replaced by a compact Riemannian manifold $\mathcal M$ of dimension $d-1$.
The results are then expressed in terms of the Ricci curvature of the manifold. Again the reader
may consult~\cite{DEL} for details.

\medskip\noindent{\bf Acknowledgment:} This work has been partially supported by the Projects STAB and Kibord
(J.D.) of the French National Research Agency (ANR). M.L. has been partially supported by
NSF Grant DMS-1301555 and the Humboldt Foundation.
\\\noindent{\scriptsize\copyright\,2016 by the authors. This paper may be reproduced, in its entirety, for non-commercial purposes.}

\bibliographystyle{siam}\bibliography{DELppn}
\end{document}